\documentclass{elsart}

\usepackage{amsmath,amssymb}
\usepackage{lscape}

\newtheorem{defi}{Definition}
\newtheorem{teo}{Theorem}

\def\adots{\mathinner{\mkern2mu\raise1pt\hbox{.}\mkern2mu\raise4pt\hbox{.}\mkern2mu\raise7pt\hbox{.}\mkern1mu}}
\begin{document}

\begin{frontmatter}
\title{Darboux tranformation and solutions of some integrable systems}
\author{D. Barrios Rolan\'{\i}a\thanksref{do}}
\ead{dolores.barrios.rolania@upm.es}
\thanks[do]{This work was partially supported by Direcci\'{o}n General de
Investigaci\'{o}n Cient\'{\i}fica y T\'ecnica, Ministerio de Econom\'{\i}a y Competitividad, under grant MTM2014-54053.}
\address{DICHOT, ETS Ingenier\'{\i}a Civil\\
Universidad Polit\'ecnica de Madrid \\
28660 Boadilla
del Monte\\ Spain}

\begin{abstract}
The relation between the Darboux transformation and the solutions of the full Kostant Toda lattice is analyzed. The discrete Korteweg de Vries equation is used to obtain such solutions and the main result of \cite{UltiAmilcar} is extended to the case of general $(p+2)$-banded Hessenberg matrices.
\end{abstract}


\end{frontmatter}

\section{Introduction}

In \cite{UltiAmilcar} some aspects of the relation between the $(p+2)$-banded matrices
\begin{equation} \label{JGeneral}
J = \left( \begin{array}{cccccc}
 a_{0,0} & 1 & &    \\
 a_{1,0} & a_{1,1} & 1 & \\
 \vdots      & \vdots    &  \ddots   & \ddots      \\
 a_{p,0} &a_{p,1} & \cdots & a_{p,p} &1 & \\
0& a_{p+1,1} & & \ddots& \ddots & \ddots\\
& 0 &\ddots  \\
&  &\ddots  \\
\end{array} \right)
\end{equation}
and the integrable system
\begin{equation}
\dot{a}_{i,j}  =  \left(a_{i,i}-a_{j,j}\right)a_{i,j} +a_{i+1,j}-a_{i,j-1}\,,\qquad i,j=0,1,\ldots
\label{sistema}
\end{equation}
were studied. In particular, a method for constructing solutions of this system was given in the case $p=2$. This method is based in the extension of the concept of Darboux transformation, which can be consulted in \cite{Bueno} for the classical tridiagonal case $p=1$. Due to the matrical interpretation of this method, the concept of transform of Darboux was extended in \cite{UltiAmilcar} for an arbitrary $p\in \mathbb{N}$ and a banded matrix $J$ as in (\ref{JGeneral}). However, just in \cite{Dani} the existence of such kind of generalized transformation was determined. As a consequence, now we are under the appropriate conditions to generalize the method for constructing solutions of (\ref{sistema}), given in \cite{UltiAmilcar} in the case $p=2$, to any $p\in \mathbb{N}$. This is precisely the object of this paper.

For simplicity of the reading, we recall here some concepts introduced in \cite{UltiAmilcar} and used in \cite{Dani} which will be used in our work. The system (\ref{sistema}) is usually called full Kostant Toda lattice. Here and in the sequel, the dot means differentiation with respect to $t\in \mathbb{R}$. However, in most of the cases we suppress the explicit $t$-dependence for brevity.

\begin{defi}
The infinite matrix $J$ is called a solution of (\ref{sistema}) if:

\begin{enumerate}
\item For each $j=0,1\ldots $ the entries $a_{i,j}=a_{i,j}(t)\,,\,i=j,j+1,\ldots, j+p\,,$ of $J$ are continuous functions with complex values defined in an open interval $\mathcal{I}_j$ such that
\begin{equation}
\bigcap_{j=0}^{N}\mathcal{I}_j \neq \emptyset \quad  \text{for any } N\in \mathbb{N}\,.
\label{intervalos1}
\end{equation}
\item The entries $a_{i,j}$ of $J$ verify (\ref{sistema}).
\end{enumerate}

\end{defi}

An important tool for us is the called discrete Korteweg de Vries (KdV) equation,
\begin{equation}
\label{VolterraGeneral}
\dot \gamma_n=\gamma_{n}\left(\sum_{i=1}^{p}\gamma_{n+i}-\sum_{i=1}^{p}\gamma_{n-i}\right)\,,\quad n\in \mathbb{N}\,.
\end{equation}
This system is an extension of the Volterra lattices studied in \cite{Peherstorfer} and \cite{Simon}. As in (\ref{sistema}), the matrices theory is used to analyze the KdV equations. The matrix associated with this system is
$$\Gamma = \left( \begin{array}{ccccc}
 0 & 1 & &    \\
 0 & 0 & 1 & \\
 \vdots  & &\ddots  & \ddots  \\
  0  & \vdots&& \\
\gamma_1 & 0  &  & \\
0& \gamma_2 &  \ddots \\
& \ddots &\ddots   \\
\end{array} \right)\,.
$$
Also we assume that the entries $\gamma_j=\gamma_j(t)\,,\,j\in \mathbb{N}\,,$ of $\Gamma$ are continuous functions with complex values defined in the open intervals $\mathcal{O}_j$ such that
\begin{equation}
\bigcap_{j=0}^{N}\mathcal{O}_j \neq \emptyset \,, \text{for any } N\in \mathbb{N}\,.
\label{intervalos2}
\end{equation}

\begin{defi}
\label{definition2}
The matrix $\Gamma$ is called a solution of (\ref{VolterraGeneral}) if the sequence $\{\gamma_n\}$ verify (\ref{VolterraGeneral}) and (\ref{intervalos2}).
\end{defi}

With respect to the extension of the Darboux transformation, the following definition was introduced in \cite{UltiAmilcar} and analyzed in \cite{Dani}. As usual, here and in the sequel $M_n$ denotes the finite matrix formed by the first $n$ rows and columns of the infinite matrix $M$.

\begin{defi}
Let $B=\left(b_{ij}\right)\,,\,i,j\in \mathbb{N}\,,$ be a lower Hessenberg $(p+2)$-banded matrix,
\begin{equation} \label{B}
B = \left( \begin{array}{cccccc}
 b_{0,0} & 1 & &    \\
 b_{1,0} & b_{1,1} & 1 & \\
 \vdots      & \vdots    &  \ddots   & \ddots      \\
 b_{p,0} &b_{p,1} & \cdots & b_{p,p} &1 & \\
0& b_{p+1,1} & & \ddots& \ddots & \ddots\\
& 0 &\ddots  \\
&  &\ddots  \\
\end{array} \right)
\end{equation}
such that $\det(B_n) \neq 0$ for any $n\in \mathbb{N}$. Let $L\,,\,U$ be two lower and upper triangular matrices, respectively, such that the entries in the diagonal of $L$ are $l_{ii}=1$ and $B=LU$ is the (unique) $LU$ factorization of $B$ in these conditions. Assume $L=L^{(1)}L^{(2)}\ldots L^{(p)}$, where
\begin{equation}
U = \left( \begin{array}{cccccc}
\gamma_1 & 1  &    \\
&  \gamma_{p+2} & 1  \\
&  &\gamma_{2p+3} & \ddots  \\
   &&& \ddots\\
\end{array} \right),\,
L^{(i)} = \left( \begin{array}{cccccc}
  1 & &    \\
 \gamma_{i+1} & 1 & \\
   & \gamma_{p+i+2}  & 1 \\
 &&\gamma_{2p+i+3} &\ddots\\
 &&&\ddots
\end{array} \right),\,i=1,\ldots,p\,.
\label{Darboux2}
\end{equation}
Then the matrical decomposition
\begin{equation}
L^{(1)}L^{(2)}\cdots L^{(p)}U
\label{(8)}
\end{equation}
is called a Darboux factorization of $B$. Moreover, any circular permutation
\begin{equation}
L^{(i+1)}\cdots L^{(p)}UL^{(1)}\cdots L^{(i)}\,,\quad i=1,2,\ldots, p\,,
\label{9}
\end{equation}
of (\ref{(8)}) is called a Darboux transformation of $B$. (We understand $UL^{(1)}\cdots L^{(p)}$ in (\ref{9}) when $i=p$.)
\end{defi}

The main goal of this paper is to obtain some solutions of (\ref{sistema}) using the Darboux transformations of $J-CI$ under certain conditions for $C\in \mathbb{C}$. More precisely, our main result is the following.
\begin{teo}\label{teorema1}
Let $J$ be a $(p+2)$-banded matrix with the structure given in (\ref{JGeneral}). Assume that $J$ is a solution of (\ref{sistema}) verifying $a_{p+i,i}\neq 0\,,\,i=0,1,\ldots$ and let $C\in \mathbb{C}$ be such that $\det(J_n-CI_n)\neq 0$ for any $n\in \mathbb{N}$. Then there exist $p$ solutions $J^{(1)},\ldots, J^{(p)}$ of (\ref{sistema}) such that the following relations hold.
\begin{equation}
J^{(i)}=CI+L^{(i+1)}\cdots L^{(p)}UL^{(1)}\cdots L^{(i)}\,,\quad i=0,1,2,\ldots,p
\label{7}
\end{equation}
(assuming $J^{(0)}=J$), where $L^{(i)}\,,\,i=1,2,\ldots,p\,,$ and $U$ have respectively the structure given in (\ref{Darboux2}). Moreover the entries of $U,\,L^{(1)}\,,\cdots \,,L^{(p)}$ provide the sequence $\{\gamma_n\}$, which defines a solution $\Gamma$ of (\ref{VolterraGeneral}).
\end{teo}

In the development of the paper it is convenient to have another expression of (\ref{7}), which is obtaining when each entry of the matrix $J^{(i)}$ is given in terms of the entries of matrices of $U,\,L^{(1)}\,,\cdots \,,L^{(p)}$. This is, if we write
\begin{equation} \label{iMatriz}
J^{(j)} = \left( \begin{array}{cccccc}
 a^{(j)}_{0,0} & 1 & &    \\
 a^{(j)}_{1,0} & a^{(j)}_{1,1} & 1 & \\
 \vdots      & \vdots    &  \ddots   & \ddots      \\
 a^{(j)}_{p,0} &a^{(j)}_{p,1} & \cdots & a^{(j)}_{p,p} &1 & \\
0& a^{(j)}_{p+1,1} & & \ddots& \ddots & \ddots\\
& 0 &\ddots  \\
&  &\ddots  \\
\end{array} \right)\,,\quad j=0,1,2,\ldots,p\,,
\end{equation}
then for each $j=0,1,\ldots, p$ from the products on the right hand side of (\ref{7}) we arrive to
\begin{eqnarray}
a^{(j)}_{i,i} & = & C+ \sum_{s=j+1}^{j+p+1} \gamma_{(i-1)p+i+s}\,,\label{10}\\
a^{(j)}_{i+k,i} & = & \sum_{E^{(j)}_{k}} \gamma_{(i-1)p+i_1+i}\gamma_{ip+i_2+i}\cdots \gamma_{(k+i-1)p+i_{k+1}+i}\,,\label{11}\\
&& \quad \quad \quad   i=0,1,\ldots\,,\quad k=1,2,\ldots\,,p\,. \nonumber
\end{eqnarray}
(we recall that $J^{(0)}=J$ and consequently $a^{(0)}_{s,r}=a_{s,r}$). The sum in (\ref{11}) is extended to the set of indices $E^{(j)}_{k}$ defined as $$E^{(j)}_{k}=\{(i_1,\ldots ,i_{k+1}):j+k+1\leq i_{k+1}\leq\cdots \leq i_1\leq j+p+1\}.$$ Frequently, relations (\ref{10})-(\ref{11}) (and by extension (\ref{7})) are known as {\it B\"{a}cklund transformations} associated with (\ref{sistema}) (see \cite{Simon}).

We underline the relevance of the Darboux factorization (\ref{(8)}) for obtaining the new solutions $J^{(i)},\,i=1,\ldots, p,$ given in (\ref{7}). In this paper we study this problem, given a method to arrive to (\ref{(8)}) from $J$.

In Section \ref{auxiliares} some tools for our work are presented and the main auxiliary results are introduced. In Section \ref{construction} an iterative method for obtaining the Darboux factorization (\ref{(8)}) is studied. Finally, the proof of Theorem 1 is in Section \ref{demostracion}.

\section{Auxiliary results}\label{auxiliares}

An important tool in our approach is the family of polynomial $\{P_n(z)\}=\{P_n(t,z)\}\,,\,n\in \mathbb{N}\,,$ associated with the matrix $J$. This family is defined by the following recurrence relation.
\begin{equation}\label{polinomiosp}
\mbox{\hspace{-.8cm}}
\left.
 \begin{array}{r}
\displaystyle\sum_{i=n-p}^{n-1}a_{n,i}P_{i}(z)+(a_{n,n}-z)P_{n}(z)+P_{n+1}(z) =0 ,\quad n=0,1,\ldots \\
 P_0(z) \equiv 1 \, ,\quad P_{-1}(z) =\cdots =P_{-p}(z)=0\,.
\end{array}
\right\}
 \end{equation}
If $J$ is a solution of (\ref{sistema}) then each polynomial $P_n=P_n(t,z)$ is a continuous function on $t$. Furthermore, in Lemma 2 of \cite{UltiAmilcar} was proved
\begin{equation}
\label{derivada}
\dot P_n(z)=-\sum_{i=n-p}^{n-1}a_{n,i}P_i\,,\quad n=0,1,\ldots \,.
\end{equation}
As a consequence of the above comments, for each $n=0,1,\ldots $ we have that $\dot P_n(z)$ is also a continuous function on $t$ in some open interval of $\mathbb{R}$.

The following lemma is obtained from (\ref{polinomiosp}) and (\ref{derivada}).

\begin{lem}
\label{lema1}
If $J$ is a solution of (\ref{sistema}) then we have
\begin{equation}\label{breve}
\dot P_n(z)=(a_{n,n}-z)P_{n}(z)+P_{n+1}(z)\,,\quad n=0,1,\ldots \,.
\end{equation}
\end{lem}

On the other hand, it is well known that for each matrix $J$ and $C\in \mathbb{C} $ in the conditions of Theorem \ref{teorema1} there exists the $LU$ factorization of $J-CI$. This is, there exists a banded lower triangular matrix
\begin{equation}\label{1-}
L= \left( \begin{array}{ccccccc}
1\\
 l_{1,1} & 1 & &    \\
 \vdots &\ddots &\ddots\\
 l_{p,1} & l_{p,2} & \dots & 1 && \\
0 & l_{p+1,2}  & \ddots & & \ddots \\
& 0 &\ddots &\ddots  &\ddots \\
\end{array} \right)\,,
\end{equation}
and there exists an upper triangular matrix $U=U(t)$ as in (\ref{Darboux2}) such that
\begin{equation}
J -CI=L U ,
\label{LU}
\end{equation}
being
\begin{equation}
J_n -CI_n=L_n U_n
\label{LUn}
\end{equation}
for each $n\in \mathbb{N}$ (see for instance \cite{Gantmacher}). The following auxiliary result also will be used in the proof of Theorem \ref{teorema1}.

\begin{lem}
\label{lema2}
In the above conditions, the entries $\gamma_{np+n+1}\,,n=0,1,\ldots, $ of $U$ verify (\ref{VolterraGeneral}).
\end{lem}
\noindent
Proof.- It is obvious that the recurrence relation (\ref{polinomiosp}) can be rewritten as
$$
(J_n-zI_n)\left(
\begin{array}{c}
P_0(z)\\P_1(z)\\\vdots \\P_{n-1}(z)
\end{array}
\right)=\left(
\begin{array}{c}
0\\0\\\vdots \\-P_{n}(z)
\end{array}
\right)\,.
$$
In particular, for $z=C$ from (\ref{LUn}) we have
$$
L_nU_n\left(
\begin{array}{c}
P_0(C)\\P_1(C)\\\vdots \\P_{n-1}(C)
\end{array}
\right)=\left(
\begin{array}{c}
0\\0\\\vdots \\-P_{n}(C)
\end{array}
\right)\,.
$$
Then, using the fact that $L$ is a triangular matrix whose diagonal entries are 1,
$$
U_n\left(
\begin{array}{c}
P_0(C)\\P_1(C)\\\vdots \\P_{n-1}(C)
\end{array}
\right)=L_n^{-1}\left(
\begin{array}{c}
0\\0\\\vdots \\-P_{n}(C)
\end{array}
\right)=\left(
\begin{array}{c}
0\\0\\\vdots \\-P_{n}(C)
\end{array}
\right)\,.
$$
This is, taking into account the structure of $U$ (see (\ref{Darboux2})),
\begin{equation}
\gamma_{np+n+1}=-\frac{P_{n+1}(C)}{P_n(C)}\,,\quad n=0,1,\ldots
\label{coci}
\end{equation}
Taking derivatives in (\ref{coci}) we arrive to
$$
\dot\gamma_{np+n+1}=-\frac{P_n(C)}{P_{n-1}(C)}\left(\frac{\dot P_n(C)}{P_{n}(C)}-\frac{\dot P_{n-1}(C)}{P_{n-1}(C)} \right)\,.
$$
From this and Lemma \ref{lema1},
$$
\dot\gamma_{np+n+1}=\gamma_{np+n+1}\left(a_{n,n}-a_{n-1,n-1}-\gamma_{np+n+1}+\gamma_{(n-1)p+n}\right)\,.
$$
Then since (\ref{10}) (with $j=0$) we arrive to (\ref{VolterraGeneral}). \hfill $\Box$

The next result guaranties the existence of the Darboux transformation, which is an important tool in the proof of Theorem \ref{teorema1}.

\begin{lem}[Theorem 1 in \cite{Dani}]
\label{teorema2}
Let $L$ be a lower triangular matrix as in (\ref{1-}) with complex entries,
such that $l_{p+j,j+1}\ne 0$ for each $j=0,1,\ldots $ Then there exists a set of $p(p-1)/2$ complex numbers
\begin{equation}
\left.
\begin{array}{ccccc}
\gamma_2, & \gamma_{p+3} &\cdots & \cdots &\gamma_{(p-1)p},\\
\gamma_3, & \gamma_{p+4} &\cdots & \gamma_{(p-2)p}, \\
\vdots & \vdots &\adots \\
\gamma_{p-1},& \gamma_{2p},\\
\gamma_p
\end{array}
\right\}
\label{puntos}
\end{equation}
and there exist $p$ triangular matrices $L^{(i)},\,i=1,\ldots,p\,$, as in (\ref{Darboux2})
such that
\begin{equation}
L=L^{(1)}L^{(2)}\cdots L^{(p)}\,,
\label{(b)}
\end{equation}
where $\gamma_{k(p+1)+i+1}\neq 0$ for $i=1,2,\ldots , p$ and $k=0,1,\ldots $ Moreover, the factorization (\ref{(b)}) is unique for each fixed set
of points (\ref{puntos}) .
\end{lem}

In \cite{Dani} also the following result is used. Here, the necessary conditions to obtain the set (\ref{puntos}) are explicit.

\begin{lem}[Theorem 2 in \cite{Dani}]
\label{teoremaDani}
Let us consider a $(p+1)$-banded lower triangular matrix $L$ as in (\ref{1-}) such that $l_{p+j,j+1}\ne 0$ for each $j=0,1,\ldots $
Assume $\alpha_1,\alpha_2,\ldots ,\alpha_{p-1}\in \mathbb{C}$ such that
\begin{equation}
\sum_{s=0}^{p-1}(-1)^s \alpha_{p-s}\alpha_{p-s+1}\cdots \alpha_{p-1}C_k^{(s)}\ne 0\,,\quad \text{for all } k=1,2,\ldots,
\label{condition}
\end{equation}
where
$C_1^{(s)}:=l_{p-s-1,1}$ and
\begin{equation}\label{*}
C_k^{(s)}:=\left|
\begin{array}{lcccl}
l_{p-s-1,1}&l_{p-s-1,2}&\cdots& \cdots&l_{p-s-1,k}\\
l_{p,1}&l_{p,2}&\ddots &\ddots &l_{p,k}\\
0 & l_{p+1,2}&\ddots &\ddots &l_{p+1,k}\\
\vdots & \ddots & \ddots& \ddots &\vdots\\
0& \cdots & 0 & l_{p+k-2,k-1}&l_{p+k-2,k}\\
\end{array}
\right|\,,\quad k\geq 2\,,
\end{equation}
for each $s=0,1,\ldots , p-1$. (We understand $\alpha_{p-s}\alpha_{p-s+1}\cdots \alpha_{p-1}=1$ for $s=0$ and $l_{i,j}=0$ for $j>i+1$.)
Then there exists a bi-diagonal matrix
\begin{equation}
D^{(1)} = \left( \begin{array}{cccccc}
  1 & &    \\
 \alpha_{1} & 1 & \\
   & \alpha_{2}  & 1 \\
 &&\alpha_{3} &\ddots\\
 &&&\ddots
\end{array} \right)
\label{LL}
\end{equation}
and there exists a $p$-banded lower triangular matrix
\begin{equation}\label{A}
A= \left( \begin{array}{ccccccc}
1\\
 \delta_{2,1} & 1 & &    \\
 \vdots &\ddots &\ddots\\
 \delta_{p,1} & \delta_{p,2} & \dots & 1 && \\
0 &\delta_{p+1,2}  & \ddots & & \ddots \\
& 0 &\ddots &\ddots  &\ddots \\
\end{array} \right)
\end{equation}
such that $\delta_{p+k-1,k}\neq 0\,,\, k=1,2,\ldots \,,$ and
\begin{equation}
L=D^{(1)}A\,.
\label{(a)(a)}
\end{equation}
Moreover, if the $p$ entries $\alpha_1,\alpha_2,\ldots ,\alpha_{p-1}\in \mathbb{C}$ of $D^{(1)}$ are fixed verifying (\ref{condition}), then (\ref{(a)(a)}) is the unique
factorization of $L$ in these conditions.
\end{lem}

Lemma \ref{teorema2} in \cite{Dani} was obtained as a corollary of Lema \ref{teoremaDani}. The key of this fact is the next lemma, which we use en the proof of our main result. Although the idea of the proof of Lemma \ref{lema5} is implicitly contained in \cite{Dani}, we include it here for an easier reading.

\begin{lem}
\label{lema5}
Let $s\in \{0,1,\ldots , p-2\}$ be and take $N\in \mathbb{N}$. Let us consider the triangular matrix
\begin{equation}
T^{(s)}= \left( \begin{array}{ccccccc}
1\\
 m^{(s)}_{1,1} & 1 & &    \\
 \vdots &\ddots &\ddots\\
 m^{(s)}_{p-s,1} & m^{(s)}_{p-s,2} & \dots & 1 && \\
0 &m^{(s)}_{p-s+1,2}  & \ddots & & \ddots \\
& 0 &\ddots &\ddots  &\ddots \\
\end{array} \right)\,,
\label{otro26}
\end{equation}
where $m^{(s)}_{p-s+j,j+1}\neq 0\,,\,j=0,1,\ldots$. Then there exist $p-s-1$ complex values
\begin{equation}
\label{D}
\alpha^{(s)}_{i}\neq 0\,, i=1,2,\ldots,p-s-1\,,
\end{equation}
such that
\begin{equation}
\sum_{j=0}^{p-s-1}(-1)^j \alpha^{(s)}_{p-s-j}\alpha^{(s)}_{p-s-j+1}\cdots \alpha^{(s)}_{p-s-1}R_k^{(s,j)}\ne 0
\label{DD}
\end{equation}
for each $k=1,2,\ldots,N ,$
where $R_1^{(s,r)}=m^{(s)}_{p-s-r-1,1}$ and
\begin{equation}\label{*}
R_k^{(s,r)}:=\left|
\begin{array}{lllcc}
m^{(s)}_{p-s-r-1,1}&m^{(s)}_{p-s-r-1,2}&\cdots& \cdots&m^{(s)}_{p-s-r-1,k}\\
m^{(s)}_{p-s,1}&m^{(s)}_{p-s,2}&\cdots& \cdots&m^{(s)}_{p-s,k}\\
0 & m^{(s)}_{p-s+1,2}&\ddots &\ddots &m^{(s)}_{p-s+1,k}\\
\vdots & \ddots & \ddots& \ddots &\vdots\\
0& \cdots & 0 & m^{(s)}_{p-s+k-2,k-1}&m^{(s)}_{p-s+k-2,k}\\
\end{array}
\right|\,,\quad k\geq 2\,.
\end{equation}
\end{lem}

\noindent
Proof.- As in the proof of Theorem 2 in \cite{Dani}, we consider the hyperplanes $\pi_k\,,\,k=1,2,\ldots, N,$ in $\mathbb{R}^{p-s-1}$ given by the equation of the form
\begin{equation}
R_k^{(s,1)}x_1+ R_k^{(s,2)}x_2+\cdots +R_k^{(s,p-s-1)}x_{p-s-1}=R_k^{(s,0)}
\label{AA}
\end{equation}
and also the hyperplanes $\Psi_i$ of equation $x_i=0\,,\,i=1,\ldots, p-s-1$. We define $\pi:=\bigcup_{k=1}^{N}\pi_k$ and $\Psi:=\bigcup_{i=1}^{p-s-1}\Psi_i$. Then, if $\mu$ is the Lebesgue measure in $\mathbb{R}^{p-s-1}$, it is well known that $\mu(\pi\cup \Psi)=0$ (see \cite{Taylor} for details). As a consequence, there exists a nonnumerable set of points $X\in \mathbb{R}^{p-s-1}$ such that
$
X=(x_1,\ldots ,x_{p-s-1})\notin \pi\cup\Psi.
$
We choose one of these points and we define iteratively
\begin{equation}
\alpha^{(s)}_{p-s-j}=
\left\{\begin{array}{ccl}
x_1 &,& \text{ if  }j=1\\ \\
\displaystyle\frac{(-1)^{j+1}x_j}{\alpha^{(s)}_{p-s-j+1}\alpha^{(s)}_{p-s-j+2}\ldots \alpha^{(s)}_{p-s-1}} & , & \text{ if  }j=2,\ldots ,p-s-1\,.
\end{array}
\right.
\label{BB}
\end{equation}
Note that $\alpha^{(s)}_{1},\ldots , \alpha^{(s)}_{p-s-1}$ are well defined because $X\notin \Psi$ and consequently $x_j\neq 0$ for each $j=1,\ldots, p-s-1$. Therefore $\alpha^{(s)}_{i}\ne 0\,,\,i=1,\ldots, p-s-1$.

From (\ref{AA}) and (\ref{BB}) we arrive to (\ref{DD}). $\hfill \Box$

\section{Construction of the matrices $L^{(i)},\,i=1,\ldots,p\,$}\label{construction}

In the conditions of Lemma \ref{teorema2}, from the set of data (\ref{puntos}) it is possible to build $L^{(1)},\ldots , L^{(p)}$ which are the factors in (\ref{(b)}). With this purpose we will use Table \ref{tabla}. From the Backl\"{u}nd transformation (\ref{11}) (for $j=0$) it is easy to arrive to
\begin{eqnarray}
\delta^{(i)}_{k}\gamma_{(k+i+1)p+i}& = &a_{k+i+1,i-1} - \sum_{
\widetilde{E}^{(0)}_{k+2}} \gamma_{(i-2)p+i+i_1-1}\gamma_{(i-1)p+i+i_2-1}\cdots \gamma_{(k+i)p+i+i_{k+3}-1}\nonumber\\
& &i\in \mathbb{N},\quad k=-1,0,1,\ldots , p-2\,,
\label{(a)}
\end{eqnarray}
where
\begin{equation}
\delta^{(i)}_{k}=\gamma_{(i-1)p+i}\gamma_{ip+i}\cdots \gamma_{(k+i)p+i}
\label{20}
\end{equation}
and
\begin{equation}
\widetilde{E}^{(0)}_{k+2}=\{(i_1,\ldots ,i_{k+3}):k+3\leq i_{k+3}\leq\cdots \leq i_1\leq p+1\,,i_{k+3}<p+1\}.
\label{3434}
\end{equation}
We remark that the subsequence $\{\gamma_{mp+m+1}\}\,,\,m\in \mathbb{N}\,,$ in (\ref{(a)})-(\ref{20}) is known and given in (\ref{coci}). This subsequence of $\{\gamma_n\}$ constitutes the first row of Table \ref{tabla} and defines the matrix $U$. With the goal to complete de rest of the $\gamma$'s, the main idea is to obtain $\gamma_{(k+i+1)p+i}$ for $k=-1,0,\ldots, p-2$, iteratively for each fixed $i\in \mathbb{N}$. Firstly we take $i=1$ in (\ref{(a)}). This is,
\begin{eqnarray}
\delta^{(1)}_{k}\gamma_{(k+2)p+1}& = &a_{k+2,0} - \sum_{
\widetilde{E}^{(0)}_{k+2}} \gamma_{-p+i_1}\gamma_{i_2}\gamma_{p+i_3}\cdots \gamma_{(k+1)p+i_{k+3}}\,,\nonumber\\
k&= &-1,0,1,\ldots , p-2\,,
\label{(c)}
\end{eqnarray}
and
$
\delta^{(1)}_{k}=\gamma_{1}\gamma_{p+1}\gamma_{2p+1}\cdots \gamma_{(k+1)p+1}.
$
We note that $\gamma_{-p+i_1}=0$ when $i_1\neq p+1$ in the sum of (\ref{(c)}). Therefore this sum can be rewritten extended to
$ \widetilde{\widetilde{E}}_{k+2} :=\{(i_2,\ldots ,i_{k+3}):k+3\le i_{k+3}\le\cdots \le i_2 \le p+1\,,i_{k+3}<p+1\}$. This is,
\begin{equation}
\delta^{(1)}_{k}\gamma_{(k+2)p+1}=a_{k+2,0} - \gamma_1\sum_{
\widetilde{\widetilde{E}}_{k+2}} \gamma_{i_2}\gamma_{p+i_3}\cdots \gamma_{(k+1)p+i_{k+3}}\,.
\label{(d)}
\end{equation}
For $k=-1$ we have $\delta^{(1)}_{-1}=\gamma_1$ and since (\ref{(d)}) we see
$$
\gamma_{p+1}=\frac{a_{1,0}}{\gamma_1} - \sum_{
j=2}^{p} \gamma_{j}\,.
$$
Now we can calculate $\delta^{(1)}_{0}=\gamma_1\gamma_{p+1}$. Taking $k=0$ in (\ref{(d)}),
$$
\gamma_{2p+1}=\frac{a_{2,0}}{\delta^{(1)}_{0}} -\frac{1}{\gamma_{p+1}} \sum_{
\widetilde{\widetilde{E}}_{2} } \gamma_{i_2}\gamma_{p+i_3}\,.
$$
In this way, in $p$ steps given for $k=-1,0,\ldots, p-2$ we obtain the entries
$$
\gamma_{p+1},\,\gamma_{2p+1},\, \ldots ,\gamma_{p^2+1},\,
$$
corresponding to the matrices $L^{(p)},\,L^{(p-1)},\ldots , L^{(1)}$ respectively (see (\ref{Darboux2})). These values of $\gamma$'s constitute the secondary diagonal in Table \ref{tabla}. Note that in the step $k+2$ it is possible to obtain $\gamma_{(k+2)p+1}$ in (\ref{(d)}) because the entries $\gamma_{i_2},\,\gamma_{p+i_3},\, \ldots ,\gamma_{(k+1)p+i_{k+3}}$ in the right hand side are in the upper triangular part of Table \ref{tabla}, over the secondary diagonal. In fact  $\gamma_{sp+i_{s+2}}$ is in the column $s+1$, including the secondary diagonal, for $s=0,1,\ldots , k$, because $i_{s+2}\leq p+1$. The last factor, $\gamma_{(k+1)p+i_{k+3}}$, is in the column $k+2$ but $i_{k+3}\leq p$ and this factor is not in the secondary diagonal, whose entry $\gamma_{(k+1)p+(p+1)}$ is being obtained in this step.

We will iterate the above procedure for constructing any parallel diagonal
$$
\gamma_{ip+i},\,\gamma_{(i+1)p+i},\, \ldots ,\gamma_{(p+i-1)p+i}
$$
in Table \ref{tabla}. Assume that we have constructed the parallel diagonals
\begin{equation}
\gamma_{sp+s},\,\gamma_{(s+1)p+s},\, \ldots ,\gamma_{(p+s-1)p+s}\,,\quad s=1,2,\ldots,i-1\,.
\label{(f)}
\end{equation}
Then we take $k=-1$ in (\ref{20}) and we have
$
\delta_{-1}^{(i)}=\gamma_{(i-1)p+i},
$
which is an entry of $U$. Therefore $\gamma_{ip+i}$ can be obtained from (\ref{(a)}) when $k=-1$. In general, for each $k\in \{-1,0,\ldots,p-2\}$ we can obtain $=\gamma_{(k+i+1)p+i}$ from (\ref{(a)}), because for this value of $k$ we have that $\delta_{k}^{(i)}$ is known in the previous step. As in the case $i=1$, the factors $\gamma_{(i-2)p+i+i_1-1},\,\gamma_{(i-1)p+i+i_2-1}\,,\ldots , \gamma_{(k+i)p+i+i_{k+3}-1}$ on the right hand side of (\ref{(a)}) are the entries of matrices $L^{(j)},\,j=1,\ldots,p\,$, corresponding to the upper triangular part of Table \ref{tabla}, which are known from the previous steps. In this form, the parallel diagonal $\gamma_{(k+i+1)p+i}\,,k=-1,0,\ldots, p-2\,,$ is obtained, corresponding with the value $s=i$ in (\ref{(f)}).

\begin{landscape}
\begin{table}
\begin{tabular}{l|lllllllllll}
$U$ & $ \gamma_{1}$ & $\gamma_{p+2}$ &
$\cdots$ &  $\gamma_{mp+(m+1)}$ & $\cdots$ & $\gamma_{(p-3)p+(p-2)}$ & $\gamma_{(p-2)p+(p-1)}$& $\cdots$ &   $\gamma_{(p+i-1)p+(i-1)} \cdots$\\\\
\hline\\
$L^{(1)}$&$\gamma_{2}$ & $\gamma_{p+3}$ &
$\cdots$ &  $\gamma_{mp+(m+2)}$ & $\cdots$ & $\gamma_{(p-3)p+(p-1)}$ & $\gamma_{(p-2)p+p}$&  $\cdots$  &  $\boldsymbol{\gamma_{(p+i-1)p+i}}$\\\\
$L^{(2)}$&$\gamma_{3}$ & $\gamma_{p+4}$ &
$\cdots$ &  $\gamma_{mp+(m+3)}$ & $\cdots$ & $\gamma_{(p-3)p+p}$ &  $\cdots$  &  $\boldsymbol{\gamma_{(p+i-2)p+i}}$\\\\
\vdots & \vdots &\vdots & &\vdots & \\\\
$L^{(p-m-1)}$&$\gamma_{p-m}$ & $\gamma_{2p-m+1}$ &
$\cdots$ &  $\gamma_{(m+1)p}$ & $\cdots$  &  $\boldsymbol{\gamma_{(m+i+1)p+i}}$\\\\
$\vdots$ & $\vdots$ &$\vdots$ &  \\\\
$L^{(p-2)}$&$\gamma_{p-1}$ & $\gamma_{2p}$ & $\cdots$ &  $\boldsymbol{\gamma_{(i+2)p+i}}$\\\\
$L^{(p-1)}$&$\gamma_{p}$ &  $\cdots$  &  $\boldsymbol{\gamma_{(i+1)p+i}}$\\\\
$L^{(p)}$&$\cdots$  &  $\boldsymbol{\gamma_{ip+i}}$
\end{tabular}
\caption{Constructing $L^{(1)},\ldots , L^{(p)}$}
\label{tabla}
\end{table}
\end{landscape}

\section{Proof of Theorem \ref{teorema1}}\label{demostracion}

If $J$ and $C$ verify the conditions of the statement in Theorem \ref{teorema1}, then (\ref{LU}) holds. Moreover for each $N\in \mathbb{N}$ there exists some open interval $\mathcal{I}_N\neq \emptyset\,,\, \mathcal{I}_N\subset \mathbb{R}\,,$ such that the entries $\gamma_{ip+i+1}(t)$ of $U$ (see (\ref{Darboux2})) and the entries $l_{i,j}(t)$ of $L\,,\,i,j\leq N\,,$ (see (\ref{1-})) are continuous functions on $\mathcal{I}_N$. In fact, these entries can be expressed in terms of products and sums of the entries of $J-CI$, which are continuous functions (see \cite{Isaacson} for instance). We assume a fixed $N\geq p-1$ in the sequel for convenience, and we set $t_0\in\mathcal{I}_N$.

In \cite{Dani}, Lemma \ref{teorema2} and the factorization (\ref{(b)}) were proved for a fixed matrix $L$ that could not depend on $t\in \mathbb{R}$. Also Lemma \ref{teoremaDani} was proved in \cite{Dani} for this kind of fixed matrices. Here, we need to extend these results for our matrix $L$. In other words, we want to prove that the first entries of the factors $L^{(i)},\,i=1,\ldots,p,$ in (\ref{(b)}) are defined in some open interval $\mathcal{I}$ such that $t_0\in \mathcal{I}\subset \mathbb{R}$. With this purpose we will apply Lemma \ref{lema5} successively for $s=0,1,\ldots ,p-2$.

In the first place we take $s=0$ and $T^{(0)}=L(t_0)$ in (\ref{otro26}). Then there exist $\alpha_1^{(0)},\ldots, \alpha_{p-1}^{(0)}$ verifying (\ref{DD}), this is
\begin{equation}
\sum_{j=0}^{p-1}(-1)^j \alpha^{(0)}_{p-j}\alpha^{(0)}_{p-j+1}\ldots \alpha^{(0)}_{p-1}R_k^{(0,j)}(t_0)\ne 0\,,\quad k=1,\ldots, N.
\label{otro37}
\end{equation}

For each $i\in \{1,2,\ldots, p-1\}$ we consider the following initial value problem,
\begin{equation}
\left.
\begin{array}{lcl}
\dot{\gamma}_{(i-1)p+(i+1)}(t) & = & \gamma_{(i-1)p+(i+1)}(t)\left( D_{i}^{(0)}(t)-\gamma_{(i-1)p+(i+1)}(t)\right)\\\\
\gamma_{(i-1)p+(i+1)}(t_0) & = & \alpha_i^{(0)}
\end{array}
\right\}
\label{F}
\end{equation}
being
$
D_{i}^{(0)}=a_{i,i}-a_{i-1,i-1}+\gamma_{(i-2)p+i}\,.
$
We analyze (\ref{F}) iteratively. If $i=1$ then $D_{1}^{(0)}=a_{1,1}-a_{0,0}$ and it is well known that (\ref{F}) has a unique solution $\gamma_2(t)$ in some open interval $\mathcal{I}^{(0)}_1$ containing $t_0$ (see \cite{Dou} for details). If $i=2$ then
$D_{2}^{(0)}=a_{2,2}-a_{1,1}+\gamma_2$ is defined in $\mathcal{I}^{(0)}_1$ and the solution $\gamma_{p+3}$ is defined in some open interval $\mathcal{I}^{(0)}_2$ containing $t_0$, being $\mathcal{I}^{(0)}_2 \subset \mathcal{I}^{(0)}_1$. Iterating this procedure, suppose that $\gamma_2, \gamma_{p+3},\ldots , \gamma_{(j-1)p+(j+1)}\,,\,j<p-1\,,$ are the solutions of (\ref{F}) which are continuous functions defined in some open interval $\mathcal{I}^{(0)}_j$ containing $t_0$. If $i=j+1$ then $D_{i}^{(0)}=a_{j+1,j+1}-a_{j,j}+\gamma_{(j-1)p+(j+1)}$ is a continuous function defined in such interval $\mathcal{I}^{(0)}_j$. Therefore (\ref{F}) has a solution $\gamma_{jp+(j+2)}$ in these conditions, which is defined in some interval $\mathcal{I}^{(0)}_{j+1} \subset \mathcal{I}^{(0)}_j$. Thus we have proved the existence of $p-1$ continuous functions
\begin{equation}
\gamma_{(i-1)p+(i+1)}\,,\quad i=1,2,\ldots, p-1\,,
\label{otro38}
\end{equation}
being these functions the respective solutions of (\ref{F}) in some open interval $\mathcal{I}^{(0)}$ containing $t_0$. We take $\mathcal{I}^{(0)}:=\mathcal{I}^{(0)}_{p-1}$. Note that in (\ref{otro38}) we have obtained the first row of (\ref{puntos}). These functions, defined in the interval $\mathcal{I}^{(0)}$, are the first $p-1$ entries in the subdiagonal of $L^{(1)}(t)$.

Another form to write (\ref{otro37}) is
\begin{equation}
\sum_{j=0}^{p-1}(-1)^j \gamma_{(p-j-1)p+(p-j+1)}(t_0)\cdots \gamma_{(p-2)p+p}(t_0)R_k^{(0,j)}(t_0)\ne 0\,,\quad k=1,\ldots, N.
\label{otro39}
\end{equation}
Due to the continuity of functions (\ref{otro38}) and $R_k^{(j)}$, it is possible to choose $\mathcal{I}^{(0)}$ sufficiently small such that (\ref{otro39}) is verified for each $t\in \mathcal{I}^{(0)}$. This is,
$$
\sum_{j=0}^{p-1}(-1)^j \gamma_{(p-j-1)p+(p-j+1)}(t)\cdots \gamma_{(p-2)p+p}(t)R_k^{(0,j)}(t)\ne 0\,,\quad k=1,\ldots, N,
$$
for $t\in \mathcal{I}^{(0)}$.

In these conditions we apply Lemma \ref{teoremaDani} and we see that $L(t)$ is factorized as in (\ref{(a)(a)}) for each $t\in\mathcal{I}^{(0)}$. This is,
\begin{equation}
L(t)=L^{(1)}(t)T^{(1)}(t)\,,\quad t\in \mathcal{I}^{(0)}\,,
\label{aaa}
\end{equation}
where $T^{(1)}$ is given in (\ref{otro26}) for $s=1$ verifying $m^{(1)}_{p-1+j,j+1}(t)\neq 0\,,t\in \mathcal{I}^{(0)}\,,j=0,1,\ldots,\,$ and $L^{(1)}$ has the structure given in (\ref{Darboux2}). We underline that the factorization (\ref{aaa}) should be understood in a formal sense, because we need to fix a $N\in \mathbb{N}$ for having the first entries of these matrices defined in some nonempty open interval. In other words, it is possible that the entries for the infinite matrices $L^{(1)}(t)$ and $T^{(1)}(t)$ are defined in $t=t_0$ but they are not simultaneously defined in a nonempty open interval. This remark should be applied in the sequel to matrices that depend on $t$, as in (\ref{otro42}) and (\ref{otro40}).

Our first goal is to prove the existence of the factors $L^{(i)}(t)\,,\,i=1,\ldots, p,$ verifying
\begin{equation}
L(t)=L^{(1)}(t)\cdots L^{(p)}(t)
\label{otro42}
\end{equation}
whose first entries are continuous functions in some open interval containing $t_0$. We proceed by induction. We assume
\begin{equation}
T^{(s)}(t)=L^{(s+1)}(t)T^{(s+1)}(t)\,,\quad t\in \mathcal{I}^{(s)}\,,
\label{otro40}
\end{equation}
for $s=0,1,\ldots, r\,,$ with $0\leq r<p-2$, being $\mathcal{I}^{(s)}\,,s=0,1,\ldots, r\,,$ open intervals such that
\begin{equation}
t_0\in \mathcal{I}^{(0)}\supseteq\mathcal{I}^{(1)}\supseteq\cdots \supseteq\mathcal{I}^{(r)}\,.
\label{otro43}
\end{equation}
In (\ref{otro40}) we assume $T^{(0)}=L$ and $T^{(s+1)}\,,\,s=0,1,\ldots, r\,,$ as in (\ref{otro26}) such that
\begin{equation}
m^{(s)}_{p-s+j,j+1}\neq 0\,,\quad j=0,1,\ldots
\label{otrotro40}
\end{equation}
Also we assume the matrices $L^{(s+1)}(t)\,,\,s=0,\ldots, r\,,$ with the structure given in (\ref{Darboux2}). This is,
\begin{equation}
L^{(s+1)}(t) = \left( \begin{array}{cccccc}
  1 & &    \\
 \gamma_{s+2}(t) & 1 & \\
   & \gamma_{p+s+3} (t) & 1 \\
 &&\gamma_{2p+s+4}(t) &\ddots\\
 &&&\ddots
\end{array} \right),
\label{otrotrotro40}
\end{equation}
being $\gamma_{s+2}(t)\,,\ldots \,,\gamma_{(p-s-1)p}(t) $ continuous functions in $\mathcal{I}^{(s)}$ such that
\begin{equation}
\dot{\gamma}_{(i-1)p+(i+s+1)}(t) = \gamma_{(i-1)p+(i+s+1)}(t)\left( D_{i}^{(s)}(t)-\gamma_{(i-1)p+(i+s+1)}(t)\right)\,,
\label{otro45}
\end{equation}
(see row $s+1$ in (\ref{puntos})) and
\begin{equation}
D_{i}^{(s)}=a_{i,i}-a_{i-1,i-1}-2\sum_{j=1}^s\gamma_{(i-1)p+(i+j)}+\sum_{j=1}^s\gamma_{ip+(i+j+1)}+
\sum_{j=0}^s\gamma_{(i-2)p+(i+j)}\,.
\label{otrotro45}
\end{equation}
We have proved the factorization (\ref{otro40}) in the above conditions for $s=0$ (see (\ref{aaa})). We want to prove that (\ref{otro40})-(\ref{otrotro45}) are verified for $s=r+1$ in some interval $\mathcal{I}^{(r+1)}$ such that $t_0\in \mathcal{I}^{(r+1)}\subseteq \mathcal{I}^{(r)}$. Since (\ref{otrotro40}) and Lemma \ref{lema5} (for $s=r+1$) we know that there exist $p-r-2$ complex values
$$
\alpha_i^{(r+1)}\neq 0\,,\quad i=1,2,\ldots , p-r-2\,,
$$
verifying
\begin{equation}
\sum_{j=0}^{p-r-2}(-1)^j \alpha_{p-r-j-1}^{(r+1)}\alpha_{p-r-j}^{(r+1)}\cdots \alpha_{p-r-2}^{(r+1)}R_k^{(r+1,j)}\ne 0\,,\quad k=1,\ldots, N\,.
\label{asterisco}
\end{equation}
Moreover we define
\begin{equation}
\alpha_i^{(s)}:=\gamma_{(i-1)p+i+s+1}(t_0)\,,\quad s=1,2,\ldots,r\,,\quad i=1,2,\ldots, p-s-1\,.
\label{otro48}
\end{equation}
We consider the following initial value problem for each $s=1,2,\ldots,r+1$ and $i=1,2,\ldots, p-s-1$,
\begin{equation}
\left.
\begin{array}{lcl}
\dot{\gamma}_{(i-1)p+i+s+1}(t) & = & \gamma_{(i-1)p+i+s+1}(t)\left( D_{i}^{(s)}(t)-\gamma_{(i-1)p+i+s+1}(t)\right)\\\\
\gamma_{(i-1)p+i+s+1}(t_0) & = & \alpha_i^{(s)}
\end{array}
\right\}
\label{E}
\end{equation}
where $D_{i}^{(s)}$ is given in (\ref{otrotro45}). Since (\ref{F}), (\ref{otro45}) and (\ref{otro48}), the continuous functions
$$
\gamma_{s+2}(t),\,\ldots, \gamma_{(p-s-1)p}(t)\,,\,t\in \mathcal{I}^{(s)}\,,\,s=0,\ldots, r\,,
$$
are the solutions of (\ref{E}) for $i=1,2,\ldots,p-s-1$ respectively. We study this initial value problem when $s=r+1$ taking $i=1,2,\ldots,p-r-2$. In the first place, if $i=1$ we have that
$$
D_{1}^{(r+1)}=a_{1,1}-a_{0,0}-2\sum_{j=1}^{r+1}\gamma_{j+1}+\sum_{j=1}^{r+1}\gamma_{p+j+2}
$$
is a continuous function defined in $\mathcal{I}^{(r)}$. Then there exists a solution $\gamma_{r+3}(t)$ of (\ref{E}) (with $i=1,\,s=r+1$), which is a continuous function in some open interval $\mathcal{I}_1^{(r+1)}$ containing $t_0$ such that $\mathcal{I}_1^{(r+1)}\subseteq\mathcal{I}^{(r)}$. Iterating the procedure, we suppose $\gamma_{(i-1)p+i+s+1}$ solutions of (\ref{E}) when $s=r+1$ and $i=1,2,\ldots, \widetilde{i}\,,$ being $\widetilde{i}<p-r-2$. Then the continuous functions
$$
\gamma_{(i-1)p+i+r+2}(t)\,,\quad i=1,2,\ldots, \widetilde{i}\,,
$$
are defined in some open interval $\mathcal{I}_{\widetilde{i}}^{(r+1)}$ containing $t_0$, being
$\mathcal{I}_{\widetilde{i}}^{(r+1)}\subseteq\mathcal{I}_{\widetilde{i}-1}^{(r+1)}\subseteq \cdots \subseteq \mathcal{I}_1^{(r+1)} $. We recall that also the functions
$$
\gamma_{(i-1)p+i+s+1}(t)\,,\quad i=1,2,\ldots,p-s-1\,,\quad s=0,1,\ldots, r,
$$
are determined in the previous steps. Therefore, $D_{\widetilde{i}+1}^{(r+1)}$ is a continuous function defined in some open interval $\mathcal{I}_{\widetilde{i}+1}^{(r+1)}$ such that $t_0\in\mathcal{I}_{\widetilde{i}+1}^{(r+1)}\subseteq\mathcal{I}_{\widetilde{i}}^{(r+1)}$. Then it is possible to set that (\ref{E}), for $i=\widetilde{i}+1$ and $s=r+1$, also has a solution in the above conditions. Thus for $s=r+1$ the continuous functions
\begin{equation}
\gamma_{(i-1)p+i+r+2}(t)\,,\quad i=1,2,\ldots,p-r-2\,,
\label{otro46}
\end{equation}
are defined in some open interval, namely $\mathcal{I}^{(r+1)}:=\mathcal{I}_{p-r-2}^{(r+1)}$, containing $t_0$. Note that $\mathcal{I}^{(r+1)}\subseteq \mathcal{I}^{(r)}$. Moreover, from (\ref{asterisco}) and the continuity of the functions (\ref{otro46}) we can assume $\mathcal{I}^{(r+1)}$ sufficiently small such that
\begin{eqnarray*}
\sum_{j=0}^{p-r-2}(-1)^j \gamma_{(p-r-j-2)p+p-j}(t)\gamma_{(p-r-j-1)p+p-j+1}(t)\cdots \gamma_{(p-r-2)p}(t)R_k^{(r+1,j)}(t)&\ne & 0\,,\\
k=1,\ldots, N\,,\quad t\in \mathcal{I}^{(r+1)}\,.
\end{eqnarray*}
Therefore we apply Lemma \ref{teoremaDani} (substituting $p$ by $p-r-1$) and we arrive to (\ref{otro40}) for $s=r+1$. Consequently, (\ref{otro40}) and the conditions given in (\ref{otro43})-(\ref{otrotro45}) are verified for $s=0,1,\ldots, p-2$. In particular,
$$
T^{(p-2)}(t)=L^{(p-1)}(t)T^{(p-1)}(t)\,,\quad t\in \mathcal{I}^{(p-2)}\,,
$$
where $L^{(p)}:=T^{(p-1)}$ is a by-diagonal matrix with the structure given in (\ref{Darboux2}). Then we obtain the factorization (\ref{otro42}) by applying successively (\ref{otro40}) in $t\in \mathcal{I}^{(p-2)}$ for $s=0,1,\ldots, p-2$. We recall that (\ref{otro42}) has a formal sense. This is, for each fixed $N\in \mathbb{N}$ there exists some open interval $\mathcal{I}_N$ such that the entries $l_{i,j}(t)\,,\,i,j\leq N\,,$ of $L(t)$ are obtained in terms of the corresponding products and sums of entries of
$
L^{(1)}\,,\,\ldots \,,\,L^{(p)}
$
for $t\in \mathcal{I}_N\cap \mathcal{I}^{(p-2)}$. However, it is possible to have
$\bigcap_{N\in \mathbb{N}}\mathcal{I}_N=\{t_0\}$. 

Now we analyze the derivative $\dot \gamma_n$ of the entries of $U$ and the recently defined matrices $L^{(i)}\,,\,i=1,\ldots, p$. In the first place we shall show that the matrix $\Gamma$ defined by the sequence $\{\gamma_n\}$ is a solution of (\ref{VolterraGeneral}) (see Definition \ref{definition2})).

The functions
\begin{equation}
\gamma_{(i-1)p+(i+s+1)}(t)\,,i=1,2,\ldots,p-s-1\,,\quad s=0,1,\ldots, p-2\,,
\label{otroA}
\end{equation}
defined for $t\in \mathcal{I}^{(p-2)}$ were obtained as the solutions of the initial value problem (\ref{E}). These functions correspond to the set (\ref{puntos}). Moreover the entries $a_{i,j}$ of $J$ can be obtained since (\ref{LU}) and (\ref{otro42}) in terms of the entries of the factors $U\,,\,L^{(i)}\,,\,i=1,\ldots, p\,.$ This is (\ref{10}) (taking $j=0$) for the diagonal entries of $J$. From this and (\ref{otrotro45}) we arrive immediately to (\ref{VolterraGeneral}) for the functions $\gamma_n$ given in (\ref{otroA}). On the other hand, (\ref{VolterraGeneral}) was proved for the entries of $U$ in Lemma \ref{lema2}. Thus we only need to verify (\ref{VolterraGeneral}) for the rest of the entries of $L^{(i)}\,,\,i=1,\ldots, p,$ that are not given in (\ref{otroA}).

The terms of the sequence $\{\gamma_n\}$  no given in (\ref{otroA}) are
\begin{equation}
\gamma_{(k+i+1)p+i}\,,\quad i=1,2,\ldots\,,\quad k=-1,0,\ldots, p-2\,.
\label{nuevas}
\end{equation}
The construction of these terms from the dates (\ref{otroA}) was analyzed in Section \ref{construction}. The functions  (\ref{otroA}), whose derivatives verify  (\ref{VolterraGeneral}), can be rewritten as in (\ref{nuevas}) taking $i=-p+2, -p+3,\ldots,-1,0$ (understanding $\gamma_n\equiv 0$ when $n\leq 0$). We recall that for each fixed $i\in N $ a parallel secondary diagonal of Table \ref{tabla} is obtained in (\ref{(a)}) for $k=-1,0,\ldots,p-2$. Next we look for an iterative expression like (\ref{(a)}) for the derivatives of these functions. This is, we want to use this construction for obtaining iteratively and simultaneously each function $\gamma_{(k+i+1)p+i}$ and its derivative $\dot\gamma_{(k+i+1)p+i}$ in terms of the functions and the derivatives functions obtained in the previous steps. In this way we assume that all the functions $\gamma_n$ on the right hand side of (\ref{(a)}) and their derivatives are known.

We shall prove the next expression for the derivative of the function $\delta^{(i)}_k$ given in (\ref{20}),
\begin{eqnarray}
\label{deri}
\dot\delta^{(i)}_k&=& \delta^{(i)}_k\left(
\sum_{j=0}^p\gamma_{(k+i)p+i+j}-\sum_{j=0}^p\gamma_{(i-2)p+i+j}\right),\nonumber \\
&&\qquad i=1,2,\ldots,\quad k=-1,0,\ldots, p-2\,.
\end{eqnarray}
We underline that for each $i, k$ the factors $\gamma_n$ of $\delta^{(i)}_k$ given in (\ref{20}) and its derivatives  $\dot\gamma_n$  verify (\ref{VolterraGeneral}) from the previous steps. From this fact,

\begin{eqnarray*}
\frac{\dot\delta^{(i)}_k}{\delta^{(i)}_k}&= &\sum_{r=-1}^k\frac{\dot\gamma_{(r+i)p+i}}{\gamma_{(r+i)p+i}} =\sum_{r=-1}^k\left(
\sum_{j=1}^p \gamma_{(r+i)p+i+j}-\sum_{j=1}^p \gamma_{(r+i-1)p+i+j-1}\right)\\
& = & \sum_{j=1}^p \gamma_{(k+i)p+i+j}+\sum_{r=0}^{k}\left(
\gamma_{(r+i-1)p+i+p}-\gamma_{(r+i-1)p+i}\right)-\sum_{j=1}^p \gamma_{(i-2)p+i+j-1}\,,\\
\end{eqnarray*}
which leads to (\ref{deri}).

Now, for $i\in N$ and $k\in \{-1,0,\ldots, p-2\}$ fixed, we study the derivatives of the terms on the right hand side of (\ref{(a)}). With this purpose for each $(i_1,i_2,\ldots,i_{k+3})\in \widetilde{E}_{k+2}^{(0)}$ we define
\begin{equation}
\Delta_k^{(i)}:=\gamma_{(i-2)p+i+i_{1}-1}\gamma_{(i-1)p+i+i_{2}-1}\cdots \gamma_{(i+k)p+i+i_{k+3}-1}\,,
\label{bbb}
\end{equation}
(see (\ref{3434})), where we assume that the derivatives of the functions $\gamma_n$ verify (\ref{VolterraGeneral}). From this we have
\begin{eqnarray}
\frac{\dot\Delta^{(i)}_k}{\Delta^{(i)}_k}
& = & \sum_{j=1}^p \gamma_{(k+i)p+i+i_{k+3}+j-1}-\sum_{j=1}^p \gamma_{(i-3)p+i+i_1+j-2} \nonumber\\
& + & \sum_{r=1}^{k+2} \left(\sum_{j=i+i_r}^{i+i_r+p-1} \gamma_{(r+i-3)p+j}-
\sum_{j=i+i_{r+1}-1}^{i+i_{r+1}+p-2} \gamma_{(r+i-3)p+j}\right)\,.
\label{ccc}
\end{eqnarray}

We are interested in (\ref{ccc}) when $(i_1,i_2,\ldots,i_{k+3})\in \widetilde{E}_{k+2}^{(0)}$ and $k\in \{-1,0,\ldots,p-2\}$, as in (\ref{(a)}).

Firstly, if $k\geq 0$ we have $i_r-i_{r+1}\leq p-k-2\leq p-2$ for $r=1,2,\ldots, k+2$. Then
\begin{equation}
i+i_{r+1}-1<i+i_r\leq i+i_{r+1}+p-2<i+i_{r}+p-1
\label{eee}
\end{equation}
for $(i_1,i_2,\ldots,i_{k+3})\in E_{k+2}^{(0)}$ and, in particular, for $(i_1,i_2,\ldots,i_{k+3})\in \widetilde{E}_{k+2}^{(0)}$. Then we simplify (\ref{ccc}) because
$$\sum_{j=i+i_r}^{i+i_r+p-1} \gamma_{(r+i-3)p+j}-\sum_{j=i+i_{r+1}-1}^{i+i_{r+1}+p-2} \gamma_{(r+i-3)p+j}$$
\begin{equation}
= \sum_{j=i+i_{r+1}+p-1}^{i+i_r+p-1} \gamma_{(r+i-3)p+j}-\sum_{j=i+i_{r+1}-1}^{i+i_{r}-1} \gamma_{(r+i-3)p+j}\,.
\label{ddd}
\end{equation}

Secondly, if $k=-1$ then in the above conditions we have $(i_1,i_2)\in \widetilde{E}_{1}^{(0)}$ and
\begin{equation}
2\leq i_2\leq i_1\leq p+1\,.
\label{ggg}
\end{equation}
Thus $i_1-i_2\leq p-1$. If
\begin{equation}
i_1-i_2\leq p-2
\label{fff}
\end{equation}
then (\ref{eee}) also is verified (with $r=1$) and consequently we arrive to (\ref{ddd}). Moreover, (\ref{fff}) holds when either $i_2\neq 2$ or $i_1\neq p+1$ in (\ref{ggg}). In this case, when $i_2=2$ and $i_1=p+1$, we arrive straight to (\ref{ddd}).

Therefore (\ref{ddd}) is verified for $r=1,2,\ldots, k+2\,,\,k=-1,0,\ldots, p-2$ and $(i_1,i_2,\ldots,i_{k+3})\in E_{k+2}^{(0)}$. From this and (\ref{ccc}) we have
\begin{eqnarray}
\frac{\dot\Delta^{(i)}_k}{\Delta^{(i)}_k} &=& \sum_{j=1}^p \gamma_{(k+i)p+i+i_{k+3}+j-1}-\sum_{j=1}^p \gamma_{(i-3)p+i+i_1+j-2}
 \nonumber\\
& + & \sum_{r=1}^{k+2} \left(\sum_{j=i+i_{r+1}+p-1}^{i+i_r+p-1} \gamma_{(r+i-3)p+j}-
\sum_{j=i+i_{r+1}-1}^{i+i_{r}-1} \gamma_{(r+i-3)p+j}\right)\,.
\label{kkk}
\end{eqnarray}
Moreover for each $j=0,1,\ldots , p$ we can show
$$\sum_{\widetilde{E}_{k+2}^{(j)}} \Delta^{(i)}_k\left(\sum_{s=i+i_{r}+p-1}^{i+i_{r-1}+p-1} \gamma_{(r+i-4)p+s}-
\sum_{s=i+i_{r+1}+p-1}^{i+i_{r}+p-1} \gamma_{(r+i-4)p+s}\right)=0\,,$$
\begin{equation}
k = -1,0,\ldots, p-2\,,\,r=2,\ldots,k+2\,.
\label{lll}
\end{equation}
In fact, in the first term of (\ref{lll}) for each $s=i+\tilde{j}+p-1\,,\,i_r\leq \tilde{j}\leq i_{r-1}$ we have
\begin{equation}
\Delta^{(i)}_k\gamma_{(r+i-4)p+s}=\gamma_{(i-3)p+i+\tilde{i}_0-1}\gamma_{(i-2)p+i+\tilde{i}_1-1}\cdots \gamma_{(i+k)p+i+\tilde{i}_{k+3}-1}
\label{llll}
\end{equation}
where
\begin{equation}
\tilde{i}_q=\left\{
\begin{array}{lll}
i_{q+1}+p &,& q=0,1,\ldots , r-2\\
\tilde{j}+p & , & q=r-1\\
i_q & , & q=r,r+1,\ldots ,k+3\,.
\end{array}
\right.
\label{mmm}
\end{equation}
In the second term of (\ref{lll}), for each $s=i+\tilde{j}+p-1$ and $i_{r+1}\leq \tilde{j}\leq i_{r}$ we have (\ref{llll})
for
\begin{equation}
\tilde{i}_q=\left\{
\begin{array}{lll}
i_{q+1}+p &,& q=0,1,\ldots , r-1\\
\tilde{j} & , & q=r\\
i_q & , & q=r+1,r+2,\ldots ,k+3\,.
\end{array}
\right.
\label{nnn}
\end{equation}
In both cases, (\ref{mmm}) and (\ref{nnn}), it is verified $(i_1,i_2,\ldots,i_{k+3})\in \widetilde{E}_{k+2}^{(j)}$ and
$
j+k+2\leq \tilde{i}_{k+3}\leq \cdots \leq \tilde{i}_{0}\leq j+2p+1\,,
$
being $\tilde{i}_{r-1}-\tilde{i}_{r}\geq p$ and taking $(\widetilde{i}_0,\widetilde{i}_1,\ldots,\widetilde{i}_{k+3})$ all the values in these conditions. Thus both sums coincide and (\ref{lll}) is verified.

Taking into account (\ref{lll}) and making some computations in (\ref{kkk}) we obtain
\begin{eqnarray}
\sum_{\widetilde{E}_{k+2}^{(0)}}\dot\Delta^{(i)}_k &=& \sum_{\widetilde{E}_{k+2}^{(0)}}\Delta^{(i)}_k\left(\sum_{j=1}^p \gamma_{(k+i)p+i+i_{k+3}+j-1}-\sum_{j=1}^{p} \gamma_{(i-3)p+i+i_1+j-2} \right.\nonumber\\
& - & \left.\sum_{j=i_2}^{i_1} \gamma_{(i-2)p+i+j-1} +\sum_{j=i_{k+3}}^{i_{k+2}} \gamma_{(k+i)p+i+j-1}\right)\,.
\label{ppp}
\end{eqnarray}

Due to $\widetilde{E}_{k+2}^{(0)}=E_{k+2}^{(0)}\setminus\{(\overbrace{p+1, \cdots ,p+1}^{(k+3)})\}$, since (\ref{ppp}) we have
\begin{eqnarray}
\sum_{\widetilde{E}_{k+2}^{(0)}}\dot\Delta^{(i)}_k &=& \sum_{E_{k+2}^{(0)}}\Delta^{(i)}_k\left(\sum_{j=1}^p \gamma_{(k+i)p+i+i_{k+3}+j-1}-\sum_{j=1}^{p} \gamma_{(i-3)p+i+i_1+j-2} \right.\nonumber\\
& - & \left.\sum_{j=i_2}^{i_1} \gamma_{(i-2)p+i+j-1}
+\sum_{j=i_{k+3}}^{i_{k+2}} \gamma_{(k+i)p+i+j-1}\right)\nonumber\\
&-&\delta^{(i)}_k\gamma_{(k+i+1)p+i}\left(\sum_{j=0}^p \gamma_{(k+i+1)p+i+j}-\sum_{j=1}^{p+1} \gamma_{(i-2)p+i+j-1}
\right)\,.
\label{pppp}
\end{eqnarray}

On the other hand we know
$(i_1,i_2,\ldots,i_{k+4})\in E_{k+3}^{(0)}$ if and only if
$
k+3\leq i_{k+4}-1\leq \cdots \leq i_{2}-1\leq i_{1}-1\leq p\,.
$
This is, 
$(i_2-1,i_3-1,\ldots,i_{k+4}-1)\in E_{k+2}^{(0)}$ and $i_2\leq i_1\leq p+1$. Then
\begin{equation}
\sum_{E_{k+3}^{(0)}}\Delta_{k+1}^{(i-1)}=
\sum_{E_{k+2}^{(0)}}\Delta_k^{(i)}\sum_{j=i_1}^p\gamma_{(i-3)p+i+j-1}\,.
\label{uu}
\end{equation}
Moreover,
$(i_1,i_2,\ldots,i_{k+4})\in E_{k+3}^{(0)}$ if and only if
$k+3\leq i_{k+3}\leq \cdots \leq i_{1}\leq p+1$ and $k+4\leq i_{k+4}\leq i_{k+3}$. This is,
$(i_1,\ldots,i_{k+3})\in E_{k+2}^{(0)}$ and $k+4\leq i_{k+4}\leq i_{k+3}$.
Then
\begin{equation}
\sum_{E_{k+3}^{(0)}}\Delta_{k+1}^{(i)}=
\sum_{E_{k+2}^{(0)}}\Delta_k^{(i)}\sum_{j=k+4}^{i_{k+3}}\gamma_{(k+i+1)p+i+j-1}\,.
\label{vv}
\end{equation}
Since (\ref{uu}) and (\ref{vv}), using (\ref{11}) we have $a_{k+i+2,i-1}-a_{k+i+1,i-2} =$
\begin{equation}
\sum_{E_{k+2}^{(0)}}\Delta_k^{(i)}
\left(\sum_{j=k+4}^{i_{k+3}}\gamma_{(k+i+1)p+i+j-1}-\sum_{j=i_1}^{p}\gamma_{(i-3)p+i+j-1}\right)\,.
\label{ww}
\end{equation}
Using again (\ref{10})-(\ref{11}), $\left(a_{k+i+1,k+i+1}-a_{i-1,i-1}\right)a_{k+i+1,k+i+1} =$
\begin{equation}
\sum_{E_{k+2}^{(0)}}\Delta^{(i)}_k\left(\sum_{j=1}^{p+1} \gamma_{(k+i)p+k+i+j+1}-\sum_{j=1}^{p+1} \gamma_{(i-2)p+i+j-1}
\right)\,.
\label{repe}
\end{equation}
From (\ref{ww}), (\ref{repe}) and (\ref{sistema}) we obtain $\dot a_{k+i+1,i-1}= $
\begin{equation}
\sum_{E_{k+2}^{(0)}}\Delta_k^{(i)}\left(\sum_{j=k+3}^{p+i_{k+3}}\gamma_{(k+i)p+i+j-1}-
\sum_{j=i_1}^{2p+1}\gamma_{(i-3)p+i+j-1}\right)\,.
\label{xx}
\end{equation}
Similarly to (\ref{lll}) it is easy to verify
\begin{equation}
\sum_{E^{(q)}_{k+2}} \Delta^{(i)}_{k}\left(\sum_{s=i_{k+3}}^{i_{k+2}}\gamma_{(k+i)p+i+s-1}-
\sum_{s=q+k+3}^{i_{k+3}}\gamma_{(k+i)p+i+s-1}\right)=0
\label{q}
\end{equation}
and
\begin{equation}
\sum_{E^{(q)}_{k+2}} \Delta^{(i)}_{k}\left(\sum_{s=i_{1}}^{q+p+1}\gamma_{(i-2)p+i+s-1}-
\sum_{s=i_2}^{i_1}\gamma_{(i-2)p+i+s-1}\right)=0
\label{r}
\end{equation}
for each $q=0,1,\ldots, p$. Taking derivatives in (\ref{(a)}) and using (\ref{deri}), (\ref{pppp}), (\ref{xx}) and (\ref{r}) (with $q=0$) we arrive to
$$
\dot\gamma_{(k+i+1)p+i}=\gamma_{(k+i+1)p+i}\left(\sum_{j=0}^{p}\gamma_{(k+i+1)p+i+j} -\sum_{j=0}^{p}\gamma_{(k+i+1)p+i-j}
\right)\,,
$$
which is (\ref{VolterraGeneral}) for $n=(k+i+1)p+i$ with $k, i$ in the mentioned conditions.

Our next target is to prove that the matrices given in (\ref{7}) are the solutions of (\ref{sistema}). With this purpose we show that the entries $\dot a_{q,r}^{(j)}$ of each $J^{(j)}$ verify (\ref{sistema}).
Because the sequence $\{\gamma_{n}\}$ verifies (\ref{VolterraGeneral}), if we take derivatives in (\ref{10}) and we make some more computations we arrive to $\dot a_{i,i}^{(j)} =$
$$
\sum_{j+1\leq i_2\le i_1\le j+p}\gamma_{(i-1)p+i+i_1+1} \gamma_{ip+i+i_2+1}- \sum_{j+1\leq i_2\le i_1\le j+p}\gamma_{(i-2)p+i+i_1} \gamma_{(i-1)p+i+i_2}
$$
which is $a_{i+1,i}^{(j)}-a_{i,i-1}^{(j)}$. Moreover, for $k=1,2,\ldots, p$ and using the notation of (\ref{bbb}), from (\ref{VolterraGeneral}) we see
\begin{eqnarray}
\dot a_{i+k,i}^{(j)} & = & \sum_{E_k^{(j)}}\Delta_{k-2}^{(i+1)}\left(\sum_{r=1}^p\gamma_{(i+k-1)p+i+i_{k+1}+r}+\right.\label{BBB}\\
&+&\left.\sum_{s=0}^{k-1}\sum_{r=1}^p\left(\gamma_{(i+s-1)p+i+i_{s+1}+r}-\gamma_{(i+s-1)p+i+i_{s+2}+r-1}\right)-
\sum_{r=1}^p\gamma_{(i-1)p+i+i_{1}-r}\right)\,.\nonumber
\end{eqnarray}
For $(i_1, \ldots , i_{k+1})\in E_{k}^{(j)}$ we have $j+k+1\leq i_{s+2}\leq i_{s+1}\leq j+p+1$ for $s=0,\ldots, k-1$. Then
$i_{s+2}\leq i_{s+1}\leq i_{s+2}+p-1\leq i_{s+1}+p$ and
\begin{eqnarray}
\sum_{r=1}^p\left(\gamma_{(i+s-1)p+i+i_{s+1}+r}-\gamma_{(i+s-1)p+i+i_{s+2}+r-1}\right)&=&\nonumber\\
\sum_{r=i_{s+2}}^{i_{s+1}}\gamma_{(i+s)p+i+r}&-&\sum_{r=i_{s+2}}^{i_{s+1}}\gamma_{(i+s-1)p+i+r}\,.
\label{novo}
\end{eqnarray}
Taking into account (\ref{novo}) in the right hand side of (\ref{BBB}), since (\ref{llll}) we have
\begin{eqnarray*}
\dot a_{i+k,i}^{(j)} & = &\sum_{E_k^{(j)}}\Delta_{k-2}^{(i+1)}\left( \sum_{r=1}^p\gamma_{(i+k-1)p+i+i_{k+1}+r}+ \right.\\
& + & \left. \sum_{r=i_{k+1}}^{i_k}\gamma_{(i+k-1)p+i+r}-\sum_{r=i_{2}}^{i_{1}}\gamma_{(i-1)p+i+r}
-\sum_{r=1}^{p}\gamma_{(i-1)p+i+i_1+r}\right).
\end{eqnarray*}
Then using (\ref{q})-(\ref{r}) we arrive to
\begin{equation}
\dot a_{i+k,i}^{(j)} = \sum_{E_k^{(j)}}\Delta_{k-2}^{(i+1)}\left( \sum_{r=j+k+1}^{i_{k+1}+p}\gamma_{(i+k-1)p+i+r}-\sum_{r=i_{1}-p}^{j+p+1}\gamma_{(i-1)p+i+r}   \right)\,.
\label{EEEE}
\end{equation}

On the other hand, from (\ref{10})-(\ref{11}) we have $(a_{i+k,i+k}^{(j)}-a_{i,i}^{(j)})a_{i+k,i}^{(j)}=$
\begin{equation}
\sum_{E_k^{(j)}}\Delta_{k-2}^{(i+1)}\left( \sum_{s=j+1}^{j+p+1}\gamma_{(i+k-1)p+i+k+s}-\sum_{s=j+1}^{j+p+1}\gamma_{(i-1)p+i+s} \right),
\label{FFF}
\end{equation}
\begin{equation}
a_{i+k+1,i}^{(j)}=\sum_{E_{k+1}^{(j)}}\gamma_{(i+k)p+i+i_{k+2}}\Delta_{k-2}^{(i+1)}\,.
\label{GGG}
\end{equation}
Moreover $(i_1, \ldots , i_{k+2})\in E_{k+1}^{(j)}$ if and only if $(\widetilde{i}_1, \ldots , \widetilde{i}_{k+1})\in E_{k}^{(j)}$ and $\widetilde{i}_{1}+1\leq i_1\leq j+p+1$, being $\widetilde{i}_r=i_{r+1}-1\,,\,r=1,\ldots,k+1$. Then
\begin{equation}
a_{i+k,i-1}^{(j)}=\sum_{E_{k+1}^{(j)}}\Delta^{(i)}_{k-1}
=\sum_{E_{k}^{(j)}}\Delta_{k-2}^{(i+1)}\sum_{s=i_1+1}^{j+p+1}\gamma_{(i-2)p+i+s-1}\,.
\label{I}
\end{equation}
Since (\ref{FFF}), (\ref{GGG}) and (\ref{I}) we have $(a_{i+k,i+k}^{(j)}-a_{i,i}^{(j)})a_{i+k,i}^{(j)}+a_{i+k+1,i}^{(j)}-a_{i+k,i-1}^{(j)}=$
\begin{equation}
=\sum_{E_k^{(j)}}\Delta_{k-2}^{(i+1)}\left( \sum_{r=j+k+1}^{p+i_{k+1}}\gamma_{(i+k-1)p+i+r}-\sum_{r=i_1-p}^{j+p+1}\gamma_{(i-1)p+i+r} \right)\,.
\label{J}
\end{equation}
Finally, comparing (\ref{EEEE}) and (\ref{J}) we arrive to (\ref{sistema}) for the entries of the matrices $J^{(j)}\,,\,j=1,\ldots,p$. With this, Theorem \ref{teorema1} is proved.  $\hfill \Box$



\begin{thebibliography}{9}

\bibitem{UltiAmilcar} D. Barrios Rolan\'{\i}a, A. Branquinho, A. Foulqui\'{e} Moreno, {\em On the full Kostant-Toda system and the discrete Korteweg-de Vries equations}, J. Math. Anal. Appl. 401 (2013), pp. 811-820.

\bibitem{Dani} D. Barrios Rolan\'{\i}a, D. Manrique, On the existence of Darboux transformations for banded matrices, Applied Mathematics and Computation 253 (2015), pp. 116-125.

\bibitem{Bueno} M.I. Bueno, F. Marcell\'an, Darboux transformation and perturbation of linear functionals, Linear Algebra and its Applications 384 (2004), pp. 215-242.

\bibitem{Dou} A. Dou, Ecuaciones Diferenciales Ordinarias, Ed. Dossat, Madrid, 1969.

\bibitem{Gantmacher} F.R. Gantmacher, The Theory of Matrices, Vol. I, AMS Chelsea Pub., Providence, Rhode Island, 2000.

\bibitem{Isaacson} E. Isaacson, H. Bishop Keller, Analysis of Numerical Methods, John Wiley \& Sons, Inc., New York.


\bibitem{Peherstorfer} F. Peherstorfer,{\em On Toda lattices and orthogonal polynomials}, J. Comput. Appl. Math. 133 (2001), 519-534.

\bibitem{Taylor} A. E. Taylor, General Theory of Functions ans Integration, Dover Pub., Inc., New York, 1985.

\bibitem{Simon} F. Gesztesy, H. Holden, B. Simon, Z. Zhao, {\em On the Toda and Kac-van Moerbeke systems}, Trans. Am. Math. Soc. 339(2) (1993), pp. 849-868.


\end{thebibliography}
\end{document}